\magnification 1200
\input amssym.def
\input amssym.tex
\parindent = 40 pt
\parskip = 12 pt

\font \medheading =cmbx7 at 16 true pt
\font \heading = cmbx10 at 13  true pt
 at 12 true pt

\def \R{{\bf R}}

\def \Z{{\bf Z}}

\centerline{\medheading Singular Integral Operators with Kernels  Associated to}
\centerline{\medheading Negative Powers of Real-Analytic Functions}
\rm
\line{}
\line{}
\centerline{\bf Michael Greenblatt}
\line{}
\centerline{May 29, 2015}
\baselineskip = 12 pt

\line{}
\line{}
\line{} 
\noindent{\heading 1. Introduction and theorems in the multiplicity one case.}

\line{}
\vfootnote{}{This research was supported in part by NSF grant  DMS-1001070} Let $n \geq 2$ and let $b(x)$ be a real-analytic function on a  neighborhood of the origin in $\R^n$ with $b(0) = 0$. By resolution of singularities, there is a number $\delta_0 > 0$ such that
on any sufficiently small neighborhood $U$ of the origin, $\int_U |f|^{-\delta} = \infty$ for $\delta \geq \delta_0$, and $\int_U |f|^{-\delta} < \infty$ for $\delta < \delta_0$. The number $\delta_0$ is sometimes referred to as the "critical integrability exponent" of $f$ at the origin. In this paper, we consider operators of the form
$$Tf(x) = \int_{\R^n} f(x - y) \,\alpha(x,y) \,m(y)\, |b(y)|^{-\delta_0} \, dy \eqno (1.1)$$
Here $\alpha(x,y)$ is a Schwartz function, and $m(y)$ is a bounded real-valued function on a neighborhood of the origin such that $m(y) |b(y)|^{-\delta_0}$ satisfies
natural derivative and cancellation conditions deriving from $b(y)$ that allows $T$ to be considered as a type of singular 
integral operator. The focus of this
paper will be to determine the boundedness properties of such $T$ on $L^p$ spaces for $1 < p < \infty$. Most of our 
results will concern the $L^2$ situation. As we will see, the operators we will consider will generalize local singular integral operators such as local versions of
Riesz transforms, and also classes of local multiparameter singular integrals.

We will see that some of our proofs immediately extend to analogues of singular Radon transforms for such singular integral
 operators. Namely, our results will cover some operators of the following form, where $x \in \R^m$ and $h$ is a real-analytic 
map from a neighborhood of the origin in $\R^n$ into $\R^m$ with $h(0) = 0$.
$$T'f(x) =  \int_{\R^n} f(x - h(y)) \,\alpha(x,y) \,m(y) \,|b(y)|^{-\delta_0} \, dy \eqno (1.2)$$

To help define what types of kernels we allow, we now delve into the resolution of singularities near the origin of a 
real-analytic function $b(x)$ with $b(0) = 0$. For this we use the resolution of
 singularities theorem of [G1], but other resolution of singularities theorems including Hironaka's famous work [H1]-[H2] can be
 used in similar ways.

By [G1], there is a neighborhood $U$ of the origin such that there exist finitely many 
coordinate change maps $\{\beta_i(x)\}_{i=1}^M$ and finitely many vectors $\{(m_{i1},...,m_{in})\}_{i=1}^M$ of nonnegative
 integers such that if $\rho(x)$ is a nonnegative smooth bump function supported in $U$ with $\rho(0) \neq 0$, then
$\rho(x)$ can be written in the form
$\rho(x) = \sum_{i=1}^M \rho_i(x)$ in such a way that each $\rho_i \circ \beta_i(x)$, after an adjustment on a set of measure
 zero, is a smooth nonnegative bump function  on a neighborhood of the origin with $\rho_i \circ \beta_i(0) \neq 0$. 
The components of each $\beta_i(x)$ are real-analytic. In addition, $\beta_i$ is a bjiection from $\{x: \rho_i \circ \beta_i(x) \neq 0, x_i \neq 0$ for all $i\}$ to $\{x: \rho_i(x) \neq 0\} - Z_i$
where $Z_i$ has measure zero, and on a connected neighborhood $U_i$ of
 the support of $\rho_i \circ \beta_i(x)$ the function $b \circ \beta_i(x)$ is well-defined and "comparable" to the monomial $x_1^{m_{i1}}...\,x_n^{m_{in}}$, meaning that there is a real-analytic function $c_i(x)$ with $|c_i(x)| > \epsilon > 0$ on 
$U_i$ such that $b \circ \beta_i(x) = c_i(x)x_1^{m_{i1}}...\,x_n^{m_{in}}$ on $U_i$.
This decomposition is such that the Jacobian determinant of $\beta_i(x)$ can be written in an analogous form
$d_i(x) x_1^{e_{i1}}...\,x_n^{e_{in}}$ on $U_i$; again the $e_{ij}$ are integers and $|d_i(x)| > \epsilon > 0$ on 
$U_i$.

\noindent  In view of the above, one has
$$\int_{\R^n} |b(x)|^{-\delta} \rho(x)\,dx = \sum_{i = 1}^M \int_{\R^n} |b(x)|^{-\delta}\rho_i(x)\,dx$$
$$= \sum_{i = 1}^M \int_{\R^n}|b \circ \beta_i (x)|^{-\delta}(\rho_i \circ \beta_i (x) )|d_i(x) x_1^{e_{i1}}...\,x_n^{e_{in}}| \,dx \eqno (1.3)$$
$$= \sum_{i = 1}^M \int_{\R^n} |c_i(x)|^{-\delta}|d_i(x)  x_1^{-\delta m_{i1} +  e_{i1}}...\,x_n^{-\delta m_{in} + e_{in}}|( \rho_i \circ \beta_i (x))\,dx \eqno (1.4)$$
Since $\rho_i \circ \beta_i (0) \neq 0$, the $i$th term of the sum $(1.4)$ is finite if each $-\delta m_{ij} + e_{ij} > -1$; that is,
if $\delta < {e_{ij} + 1 \over m_{ij}}$.
Thus the number $\delta_0$ is given in terms of the resolution of singularities of $b(x)$ by $\inf_{i,j}  {e_{ij} + 1 \over m_{ij}}$.
The following notion plays a major role in this paper.

\noindent {\bf Definition 1.1.} The {\it multiplicity} of the critical integrability exponent $\delta_0$ of $b(x)$ at the origin is the  maximum over all $i$ of the cardinality of $\{j: {e_{ij} + 1 \over m_{ij}} = \delta_0\}$.

One example of the significance of the multiplicity is as follows. Let $B_r$ denote $\{x \in \R^n: |x| < r\}$ and let $m$ denote the multiplicity of the exponent $\delta_0$ for $b(x)$ at the origin. It can be shown (see chapter 7 of [AGV] for details) that if  $r > 0$ is sufficiently small then as $\epsilon \rightarrow 0$  one has asymptotics of the form
$$|\{x \in B_r: |b(x)| < \epsilon\}| = c_r \epsilon^{\delta_0} (\ln \epsilon)^{m-1} + o(\epsilon^{\delta_0} (\ln \epsilon)^{m-1}) \eqno (1.5)$$
Here $c_r > 0$. One obtains analogous asymptotics for various oscillatory integrals associated with $b(x)$. Note that $(1.5)$ shows
that the multiplicity is independent of the which resolution of singularities process is being used.

\noindent {\bf Singular integrals associated to negative powers of multiplicity one functions.}

Let $b(x)$ be a real-analytic function on a neighborhood of the origin, not identically zero, with $b(0) = 0$.  Let $\rho(x)$, $\{\beta_i(x)\}_{i=1}^M$,  and
$\{(m_{i1},...,m_{in})\}_{i=1}^M$ be as above. We will define singular integrals associated to $b(x)$ as follows. For a given $i$,
 we move into the "blown-up" coordinates determined by $\beta_i(x)$ and define a type of singular integral
that is of magnitude bounded by $ C|x_1^{m_{i1}}...\,x_n^{m_{in}}|^{-\delta_0}$, with corresponding bounds on first derivatives, which is supported on the 
support of $\rho_i \circ \beta_i (x)$. An appropriate cancellation condition will be assumed that will ensure that the kernels
are distributions. A singular integral associated to $b(x)$ will then be defined to be a sum from $i=1$ to $i = M$ of the blow-downs of
such singular integrals into the original coordinates.

Specifically, we consider $k_i(x) = \sum_{(j_1,...,j_n) \in \Z^n} k_{i, j_1,...,j_n}(x)$, where for some fixed $C_0$
the function $k_{i, j_1,...,j_n}(x)$ is $C^1$, supported on $\{x: |x_l| \in [2^{-j_l}, C_02^{-j_l}]$ for all $l\}$, and satisfies
$$|k_{i, j_1,...,j_n}(x_1,...,x_n)| < C_0 | x_1^{m_{i1}}...\,x_n^{m_{in}}|^{-\delta_0} \eqno (1.6)$$
We also assume that for each $l = 1,...,n$ we have
$$|\partial_{x_l} k_{i, j_1,...,j_n}(x_1,...,x_n)| <  C_0 {1 \over |x_l|}|  x_1^{m_{i1}}...\,x_n^{m_{in}}|^{-\delta_0} \eqno (1.7)$$
The cancellation condition we assume for the multiplicity one case is that for some $\epsilon_0 > 0$, whenever $i$ and $l$ are such
that ${e_{il} + 1 \over m_{il}} = \delta_0$ (the minimum possible value), then where $Jac_{\beta_i}(x)$ denotes the Jacobian determinant of $\beta_i$ we have
$$\bigg|\int_{\R} k_{i, j_1,...,j_n}(x_1,...,x_n) \,Jac_{\beta_i}(x_1,...\,x_n)\, dx_l\bigg| < C_0 2^{- \epsilon_0 j_l} \eqno (1.8)$$
To ensure that our singular integrals are well-defined, we also 
assume that the support of $k_{i, j_1,...,j_n}(x)$ is contained in that of $\rho_i \circ \beta_i (x)$. We next make the following definition.

\noindent {\bf Definition 1.2.} If $b(x)$ has multiplicity one at the origin, we define a {\it singular integral kernel associated to $b(x)$} to
be a function $K(x)$ of the form $K(x) = \sum_{i = 1}^M \rho_i(x) k_i(\beta_i^{-1}(x))$, where $k_i$ satisfies $(1.6)-(1.8)$ and the support condition stated afterwards.

One can simply explicitly construct $K(x)$ satisfying Definition 1.2 for any given $b(x)$ with multiplicity one at the origin, but a familiar example can be derived
from local Riesz transforms:

\noindent {\bf Example.} Let $L(x)$ be the local Riesz transform kernel given by $\phi(x) {x_l \over |x|^{n+1}}$ for a cutoff 
function $\phi(x)$ supported near the origin. Then
$L(x)$ satisfies Definition 1.2 for $b(x) = x_1^2 + ... + x_n^2$. Here $\delta_ 0 = {n \over 2}$. For one can write 
$L(x) = \sum_{i = 1}^n L_i(x)$ where
$L_i(x)$ is supported on a cone centered at the $x_i$-axis. Then if $\beta_i(x) = (x_ix_1,...,x_ix_{i-1},x_i, x_ix_{i+1},...,x_ix_n)$,
the functions $L_i \circ \beta_i(x)$ will satisfy $(1.6)-(1.8)$ with $ x_1^{m_{i1}}...\,x_n^{m_{in}} = x_i^{-n}$. Nonisotropic
versions of $L(x)$ will satisfy Definition 1.2 for $b(x)$ of the form $x_1^{2k_1} + ... + x_n^{2k_n}$ for positive integers 
$k_1,...,k_n$.

Each $K(x)$ satisfying Definition 1.2 can be viewed in a natural way as a distribution as follows. Let $k_{iL}(x)$ denote the truncated version of
$k_i(x)$ given by 
$$k_{iL}(x) = \sum_{(j_1,...,j_n) \in \Z^n:\, j_l < L \,\,for\,\,all\,\, l} k_{i, j_1,...,j_n}(x) \eqno (1.9)$$
Define the corresponding truncated $K_L(x)$ by  $K_L(x) = \sum_{i = 1}^MK_{iL}(x)$, where $K_{iL}(x) =  \rho_i(x) k_{iL}(\beta_i^{-1}(x))$. Note that $K_L(x)$ is
a smooth compactly supported function. If $\phi(x)$ is a Schwartz function, then one has
$$\int_{\R^n} K_{L}(x)\,\phi(x)\,dx = \sum_{i = 1}^M \int_{\R^n}  \rho_i(x) \,k_{iL}(\beta_i^{-1}(x))\,\phi(x)\,dx$$
$$=  \sum_{i = 1}^M \int_{\R^n}  (\rho_i \circ \beta_i (x)) \,k_{iL}(x)\,Jac_{\beta_i}(x)\,(\phi \circ \beta_i (x))\,dx \eqno (1.10)$$
$$= \sum_{i = 1}^M \int_{\R^n} \sum_{ (j_1,...,j_n): j_l < L \,\,for\,\,all\,\,  l} ( \rho_i \circ \beta_i (x) )\,k_{i,j_1,...,j_n}(x)\, Jac_{\beta_i}(x)(\phi \circ \beta_i (x))\,dx \eqno (1.11)$$
Let $\psi_i(x) = (\rho_i \circ \beta_i (x))(\phi \circ \beta_i (x))$. Then $\psi_i(x)$ is a smooth compactly supported function. Note
that we do use the fact from [G1] that the function $\beta_i(x)$ is defined and smooth on a neighborhood of the support of 
$\rho_i \circ \beta_i (x)$ so that there are no issues concerning the smoothness of  $\phi \circ \beta_i (x)$ on the boundary of the support of $\rho_i \circ \beta_i (x)$. Thus we can rewrite the expression $(1.11)$ for $\int_{\R^n} K_L(x)\phi(x)\,dx$ as
$$= \sum_{i = 1}^M \int_{\R^n}\sum_{ (j_1,...,j_n): j_l < L \,\,for\,\,all\,\,  l} k_{i,j_1,...,j_n}(x)\, Jac_{\beta_i}(x)\, \psi_i(x)\,dx \eqno (1.12)$$
If $i$ is such that $|x_1^{m_{i1}}...\,x_n^{m_{in}}|^{-\delta_0}|Jac_{\beta_i}(x)| \sim|x_1^{m_{i1}}...\,x_n^{m_{in}}|^{-\delta_0} |x_1^{e_{i1}}...\,x_n^{e_{in}}|$ is integrable on a neighborhood of the origin, then by $(1.6)$, the form of the $i$th term of $(1.12)$ ensures that the kernel $K_{iL}(x)$ is a distribution that converges as $L \rightarrow \infty$ to a finite measure which we denote by $K_i(x)$.

Next, we show that for the $i$ for which  $|x_1^{m_1}...\,x_n^{m_n}|^{-\delta_0}|Jac_{\beta_i}(x)|$ is not integrable, the
cancellation condition $(1.8)$ ensures that such an $K_{iL}$ too converges as $L$ goes to infinity in the distribution sense to 
some $K_i(x)$. We will then define $K(x) = \sum_{i=1}^M K_i(x)$. To see why this is the case, note that since $b(x)$ has multiplicity one, for each
such $i$ there is exactly one value $l_0$ for which ${e_{il_0} + 1 \over m_{il_0}} = \delta_0$, and ${e_{il} + 1 \over m_{il}} > \delta_0$ for all other values of $l$. Write $\psi_i(x) = \psi_i(x_1,...\,x_{l_0 - 1}, 0, x_{l_0 + 1},...,x_n) + 
x_{l_0} \xi_i(x_1,...,x_n)$, with $\xi_i$ smooth.

The $i$th term of $(1.12)$ can be written as the sum of two terms. In the 
first, $\psi_i(x)$ is replaced by  $\psi_i(x_1,...\,x_{l_0 - 1}, 0, x_{l_0 + 1},...,x_n)$ and in the second $\psi_i(x)$ is
replaced by $\xi_i(x)$ and $k_{i,j_1,...,j_n}(x)$ by $x_{l_0}k_{i,j_1,...,j_n}(x)$. The second term is handled exactly as 
we handled the terms  for which  $|x_1^{m_1}...\,x_n^{m_n}|^{-\delta_0}|Jac_{\beta_i}(x)|$ is integrable since the 
additional $x_{l_0}$ factor causes us to once again have absolute integrability of the limiting kernel. As for the first term, we 
perform the $x_{l_0}$ integration first in the $i$th term of $(1.12)$. The cancellation condition $(1.8)$ implies that the
limiting kernel of the result of this integration is similarly absolutely integrable in the remaining $n-1$ variables, and thus the
limit again defines a distribution. 

Thus we see that $K_i(x)$ is a well defined distribution for all $i$ and therefore $K(x) = \sum_{i=1}^M K_i(x)$ gives a well-defined 
distribution. Hence if $\alpha(x,y)$ is a Schwartz function on $\R^{n + m}$ and $f(x)$ is a Schwartz function on $\R^n$ then $Tf (x)= \int_{\R^n} f(x - y) \alpha(x,y) K(y)\,dy$ is well-defined. If for some $1 < p < \infty$ and some constant $C$ the operators $T_Lf (x)= \int_{\R^n} f(x - y) \alpha(x,y) K_L(y)\,dy$ are such that $||T_L||_{L^p \rightarrow L^p} \leq C$ for all Schwartz functions $f$ and all $L$, then an application of the dominated convergence theorem gives that  one also  has $||Tf||_{L^p(\R^n)} \leq C||f||_{L^p(\R^n)}$ for all Schwartz functions, for the same  constant $C$.

\noindent Our first theorem is simply that $T$ is bounded on $L^2(\R^n)$.

\noindent {\bf Theorem 1.1.} Whenever the critical integrability exponent of $b(x)$ at the origin has multiplicity one, then there is a neighborhood $U$ of the origin such that if $K(x)$ is supported in $U$, there is a constant $C$ such that for all Schwartz functions $f(x)$ one has $||Tf||_{L^2(\R^n)} \leq C||f||_{L^2(\R^n)}$.

It turns out that it is no harder prove $L^2$ boundedness for singular Radon transform 
generalizations of $T$. Namely, let $K(x)$ be as above, and let $h_1(x),...,h_m(x)$ be real-analytic functions on a neighborhood 
of the origin in $\R^n$ with $h_i(x) = 0$ for all $i$. Let $\alpha(x,y)$ be a Schwartz function on $\R^m \times \R^n$. Then 
for a Schwartz function $f(x)$ in $m$ variables, we define $T'f(x)$ by
$$T'f(x) = \int_{\R^n} f(x_1 - h_1(y),..., x_m - h_m(y)) \,\alpha(x,y)\,K(y)\, dy \eqno (1.13)$$
The operator $T$ above corresponds to $m = n$ and $h_i(y) = y_i$ for all $i$. We have the following theorem.

\noindent {\bf Theorem 1.2.} Whenever the critical integrability exponent of $b(x)$ at the origin has multiplicity one, then there is a neighborhood $U$ of the origin such that if $K(x)$ is supported in $U$, there is a constant $C$ such that for all Schwartz functions $f(x)$, one has $||T'f||_{L^2(\R^{m})} \leq C||f||_{L^2(\R^{m})}$

 \noindent To give a rough idea of how our proofs will work, note that $(1.13)$ can be written as
$$\sum_{i = 1}^M  \int_{\R^n} f(x_1 - h_1(y),..., x_m - h_m(y))  \, \alpha(x,y)\, \rho_i(y)\, k_i(\beta_i^{-1}(y))\, dy \eqno (1.14)$$
Let $T_i$ be the operator corresponding to the $i$th term of $(1.14)$. Doing a change of variables from $y$ to $\beta_i(y)$ in the integral of $(1.14)$ leads to
$$T_i f(x) = \int_{\R^n} f(x_1  - h_1 \circ \beta_i(y), ..., x_m  - h_m \circ
\beta_i(y))\,\alpha (x,\beta_i(y)) (\rho_i \circ \beta_i(y)) k_i(y)\, Jac_{\beta_i}(y)\,dy \eqno (1.15)$$
$T_i$ is a sort of singular Radon transform with kernel $ \alpha (x,\beta_i(y))\,( \rho_i \circ \beta_i(y))\, k_i(y) \, Jac_{\beta_i}(y)$,
which we will be able to analyze by reducing to singular Radon transform estimates the author used in [G3]. 

\noindent {\heading 2. Theorems when the multiplicity is greater than one.}

When the critical integrability exponent $\delta_0$ has multiplicity greater than one at the origin, the coordinate changes $\beta_i(x)$ used in the multiplicity one case will lead to trying to
prove $L^2$ boundedness of an operator that resembles a multiparameter singular Radon transform, rather than a (one-parameter)
singular Radon transform.  Unfortunately since the $\beta_i(x)$ here involve blowups, one often ends out with a multiparameter singular Radon transform that is not bounded on $L^2$. As a result, instead 
of trying to find a general correct notion of singular integral and prove a general result, when the multiplicity is greater than one we will focus on theorems
that can be proven in the original coordinates.

\noindent {\bf Newton polyhedra and related matters.}

One can often determine the criticial integrability exponent of a function at the origin and its mutliplicity through the use of Newton polyhedron of the function. We turn to the relevant definitions.

\noindent {\bf Definition 2.1.} Let $b(x)$ be a real-analytic function with Taylor series
$\sum_{\alpha} b_{\alpha}x^{\alpha}$ on a neighborhood of the origin.
For each $\alpha$ for which $b_{\alpha} \neq 0$, let $Q_{\alpha}$ be the octant $\{t \in \R^n: 
t_i \geq \alpha_i$ for all $i\}$. The {\it Newton polyhedron} $N(b)$ of $b(x)$ is defined to be 
the convex hull of all $Q_{\alpha}$.  

A Newton polyhedron can contain faces of various dimensions in various configurations. The faces can be either compact or unbounded. In this paper, as in earlier work such as [G2] and [V], an 
important role is played by the following functions, associated to each compact face of $N(b)$.  We consider each vertex
of $N(b)$ to be a compact face of dimension zero.

\noindent {\bf Definition 2.2.} Let $F$ be a compact face of $N(b)$. Then
if $b(x) = \sum_{\alpha} b_{\alpha}x^{\alpha}$ denotes the Taylor expansion of $b$ like above, 
we define $b_F(x) = \sum_{\alpha \in F} b_{\alpha}x^{\alpha}$.

\noindent We will also use the following terminology.

\noindent {\bf Definition 2.3.} Assume $N(b)$ is nonempty. Then the 
{\it Newton distance} $d(b)$ of $b(x)$ is defined to be $\inf \{t: (t,t,...,t,t) \in N(b)\}$.

\noindent {\bf Definition 2.4.} The {\it central face} of $N(b)$ is the face of $N(b)$ of minimal dimension intersecting the line
$t_1 = t_2 = ... = t_n$. 

\noindent In Definition 2.4, the central face of $N(b)$ is well-defined since it is the intersection of all faces of $N(b)$
intersecting the line $t_1 = t_2 = ... = t_n$. An equivalent definition that can be used (such as in [AGV]) is that the central
face of $N(b)$ is the unique face of $N(b)$ that intersects the line $t_1 = t_2 = ... = t_n$ in its interior.

Extending results of [V], in [G2] the author showed that if the zeros of each $b_F(x)$ on $(\R - \{0\})^n$ are of order less than $d(b)$, then the critical integrability index $\delta_0$ is equal to ${1 \over d(b)}$ and the multiplicity is equal to $n$ minus the dimension of the central face of $N(b)$. This can be used to compute $\delta_0$ and its multiplicity for specific examples of interest, such as in the following two examples (which are covered by [V]).
Suppose $b(x) = x_1^{k_1} + ... + x_n^{k_n}$ with each $k_i$ even. Then $\delta_0 = ({1 \over k_1} + ... + {1 \over k_n})^{-1}$ and $m = 1$. On the other hand, if $b(x) = x_1^{l_1}...\,x_n^{l_n}$ then $\delta_0 = {1 \over \max_i l_i}$ and $m$ is equal to the number of times $k$ that $\max_i l_i$ appears in $\{l_1,...,l_n\}$. For in the former case the line $t_1 = ... = t_n$ intersects $N(b)$ in the
interior of the $n-1$ dimensional face with equation ${t_1 \over k_1} + ... + {t_n \over k_n} = 1$, while in the latter case the line
$t_1 = ... = t_n$ intersects $N(b)$ in the $n - k$ dimensional plane determined by the equations $t_l = \max_i l_i$ for all $l$
such that $l = \max_i l_i$. 

\noindent {\bf The function $b^*(x)$}

In order to understand the behavior of functions satisfying the finite-type condition of [G2], it is often helpful to consider the
function  $b^*(x)$ defined by
$$b^*(x) = \sum_{(v_1,...,v_n)\,\,a \,\,vertex \,\,of\,\,N(b)} |x_1^{v_1}...\,x_n^{v_n}| \eqno (2.1)$$
By Lemma 2.1 of [G2], there
is a constant $C$ such that for all $x$ one has $|b(x)| \leq Cb^*(x)$. In Lemma 4.1 of this paper we will see that
given any $\delta > 0$ there is a $\delta' > 0$ such that $|b(x)| > \delta' b^*(x)$ on a portion of any dyadic rectangle with
measure at least $1 - \delta$ times that of the rectangle. Hence $|b(x)| \sim b^*(x)$ except near the zeroes of
$|b(x)|$. 

Next, observe that the Newton polygon of any first partial $\partial_{x_l} b(x)$ is a subset of the shift of $N(b)$ by $-1$ units in 
the $x_l$ direction. Hence the above considerations tell us that
$$|\partial_{x_l} b(x)| \leq C{1 \over |x_l|} b^*(x) \eqno (2.2a)$$
If $b(x) \neq 0$, we also have
$$|\partial_{x_l} (|b(x)|^{-\delta_0})| \leq C{1 \over |x_l|} b^*(x)|b(x)|^{-1-\delta_0} \eqno (2.2b)$$

\noindent {\bf Singular integrals when the multiplicity is greater than one.}

 When the multiplicity is greater than one, the class of $b(x)$ where we will prove $L^2$ boundedness of 
associated singular integrals are the $b(x)$ analyzed in [G2] that were discussed above $(2.1)$. Namely, using the terminology of Definitions 2.1-2.4, 
we will assume that for each compact face $F$ of $N(b)$, each zero of  each $b_F(x)$ in $(\R - \{0\})^n$ has order less than 
$d(b)$. As mentioned above, by $[G2]$ in this situation we have $\delta_0 = {1 \over d(b)}$. Note that our theorems 
do not require the multiplicity to be greater than one, and in fact the theorems here
 will include some multiplicity one operators not covered by Theorem 1.1.

In the situation at hand, we define a singular integral associated to $b(x)$ as follows. 
 Let $\alpha(x,y)$ be a Schwartz function on $\R^{n + n}$. We
will consider kernels of the form $\alpha(x,y)K(y)$, where $K(y)$ is as follows. We assume that $K(y)$ can be written as 
$\sum_{(j_1,...,j_n) \in \Z^n}  K_{j_1,...,j_n}(y)$ such that for some fixed $C_1$
the function $K_{j_1,...,j_n}(y)$ is supported on  $\{x: |x_l| \in [2^{-j_l}, C_12^{-j_l}]$ for all $l\}$, is $C^1$ on $\{x: b(x) \neq 0\}$, and satisfies the estimates
$$|K_{j_1,...,j_n}(y_1,...,y_n)|  <  C_1|b(y_1,...,y_n)|^{-\delta_0} \eqno (2.3)$$
Motivated by $(2.2b)$, we also assume that if $b(y_1,...,y_n) \neq 0$ then for each $l$ we have
$$|\partial_{y_l}K_{j_1,...,j_n}(y_1,...,y_n)| <  C_1 {1 \over |y_l|} b^*(y_1,...,y_n)|b(y_1,...,y_n)|^{-1 - \delta_0} \eqno (2.4)$$

\noindent We further assume the cancellation conditions that for each $l$ we have
$$\int_{\R} K_{j_1,...,j_n}(y_1,...,y_n)\,dy_l = 0 \eqno (2.5)$$
In Lemma 4.2,
we will see that in the settings of our theorems (Theorems 2.1 and 2.2) each $K_{j_1,...,j_n}(y_1,...,y_n)$ is integrable, so
$(2.5)$ makes sense if we assume it holds whenever $K_{j_1,...,j_n}(y_1,...,y_n)$ is integrable in the $y_l$ variable for fixed
values of the other $y$ variables.

\noindent We will also assume that $K_{j_1,...,j_n}(y_1,...,y_n)$ is identically zero when $[2^{-j_1}, C_12^{-j_1}] \times$ ... $\times [2^{-j_n}, C_12^{-j_n}]$ is not contained in a certain neighborhood of the origin to be determined by our arguments. 

Some
motivation for our definition of a singular integral associated to $b(x)$ is the fact that 
for traditional multiparameter singular integrals, often a sufficient and necessary condition for $L^p$ boundedness is that 
the kernel be expressible as a dyadic sum of terms satisfying standardized estimates as well as a cancellation condition. We refer
to Theorem 2 of Lecture 2 of [N] for an example of a theorem of this nature. We also refer to the standard references [FS] and [NS]
for more general information about multiparameter singular integrals.

\noindent {\bf Example 1.} Let $b(x) = x_1^{a_1}...\,x_n^{a_n}$ for nonnegative integers $a_1,...,a_n$ with at least one $a_i$ being nonzero. Then $\delta_0 = {1 \over \max_i a_i}$ here. If $\phi(x)$ is a cutoff function supported on a sufficiently small neighborhood of the origin, $K(x) = (-1)^{sgn(x_1) + ... + sgn(x_n)} \phi(x_1^2,...,x_n^2) |b(x)|^{-\delta_0}$ will satisfy
$(2.3)-(2.5)$. In particular, $K(x) =  (-1)^{sgn(x_1) + ... + sgn(x_n)} \phi(x_1^2,...,x_n^2) {1 \over |x_1...\,x_n|}$ 
qualifies.

\noindent {\bf Example 2.} Let $f(x_1,...,x_n)$ be any real-analytic function with $f(0,...,0) = 0$, and let $b(x) = f(x_1^2,...,x_n^2)$. Then
if $\phi(x)$ is a cutoff function supported on a sufficiently small neighborhood of the origin, $K(x) = (-1)^{sgn(x_1) + ... + sgn(x_n)} \phi(x_1^2,...,x_n^2) |b(x)|^{-\delta_0}$ will satisfy  $(2.3)-(2.5)$.

For a fixed value of $x$, the function $\alpha(x,y)K(y)$ can be viewed in a natural way as a distribution in the $y$ variable as follows. This will resemble the discussion following $(1.9)$. Let $K_L(y) = \sum_{j_l < L\,\,for\,\, all\,\, l} K_{j_1,...,j_n}(y)$ and let $\phi(y)$ be a Schwartz function. Then 
$$\int_{\R^n} \alpha(x,y)\,K_L(y)\,\phi(y)\,dy = \int_{\R^n} \sum_{j_l <  L\,\,for\,\, all\,\, l}   K_{j_1,...,j_n}(y)\,\alpha(x,y)\,\phi(y)\,dy_1...dy_n \eqno (2.6)$$
Let $\sigma_x(y) = \alpha(x,y)\,\phi(y)$. Then we may write $\sigma_x(y)  = \sigma_x(0,y_2,...,y_n) + y_1 \xi_x(y_1,...,y_n) $
for some smooth $\xi_x(y_1,...,y_n)$. Then the right-hand side of $(2.6)$ can be rewritten as
$$ \int_{\R^n} \sum_{j_l < L\,\,for\,\, all\,\, l}   K_{j_1,...,j_n}(y_1,...,y_n)\,\sigma_x(0,y_2,...,y_n) \,dy_1...dy_n $$
$$+ \int_{\R^n} \sum_{j_l < L\,\,for\,\, all\,\, l}   K_{j_1,...,j_n}(y_1,...,y_n)\,y_1\, \xi_x(y_1,...,y_n) \,dy_1...dy_n \eqno (2.7)$$
Because of the cancellation condition $(2.5)$ in the $y_1$ variable, the first integral of $(2.7)$ is zero. We next similarly write
$\xi_x(y_1,...,y_n) = \xi_x(y_1,0,y_3,...,y_n) + y_2 \tilde{\xi}_x(y_1,...,y_n)$ and insert it into $(2.7)$, obtaining
$$ \int_{\R^n} \sum_{j_l < L\,\,for\,\, all\,\, l}   K_{j_1,...,j_n}(y_1,...,y_n)\,y_1y_2 \,\tilde{\xi}_x(y_1,...,y_n) \,dy_1...dy_n \eqno (2.8)$$
Going through all the $y_l$ variables in this way, we see that $\int_{\R^n} \alpha(x,y)\,K_L(y)\,\phi(y)\,dy$ is equal to an
expression 
$$ \int_{\R^n} \sum_{j_l < L\,\,for\,\, all\,\, l}   K_{j_1,...,j_n}(y_1,...,y_n)\,y_1y_2...y_n \,\eta_x(y_1,...,y_n) \,dy_1...dy_n \eqno (2.9)$$
Here $\eta_x(y_1,...,y_n)$ is smooth in both the $x$ and $y$ variables. We will see in Lemma 4.2 that the condition on the order of the zeroes of
the functions $b_F(x)$ on $(\R - \{0\})^n$ implies that the integral of $|b(x)|^{-\delta_0}$ over any dyadic rectangle in $U$ is uniformly bounded. Thus $(2.3)$ implies that  $K_{j_1,...,j_n}(y)\,y_1y_2...y_n$ is absolutely integrable over $U$. Hence $\alpha(x,y)K(y)$ is naturally a distribution in $y$ when $K(y)$ is 
supported in $U$ if its action on $\phi(y)$ is given by
$$\langle \alpha(x,y)K(y), \phi(y) \rangle=  \int_{\R^n} \sum_{(j_1,...,j_n) \in \Z^n}  K_{j_1,...,j_n}(y_1,...,y_n)\,y_1y_2...y_n \,\eta_x(y_1,...,y_n) \,dy_1...dy_n \eqno (2.10)$$
One can then use $(2.10)$ to define $Tf(x) = \int_{\R^n} f(x - y)\alpha(x,y)K(y)\,dy$ for Schwartz functions $f$, and then
examine boundedness of such integral operators on $L^p$ spaces. We have the following theorem in this regard for $p = 2$.

\noindent {\bf Theorem 2.1.} Suppose each polynomial $b_F(x)$ only has zeroes of order less than $d(b)$ on $(\R - \{0\})^n$. Suppose also that there is a $C_0 > 0$ and a neighborhood $U$ of the origin such that for each $l$, there is a set $Z \subset \R^{n-1}$ of measure zero such that the function $\partial_{x_l} b(x)$ has at most $C_0$ zeroes in $U$ 
on any line parallel to the $x_l$ coordinate axis whose projection onto the plane $x_l = 0$ is not in $Z$.
Then there is an $R > 0$ such that if each $K_{j_1,...,j_n}(y)$ satisfies $(2.3)-(2.5)$ and is supported on $|y| < R$, then 
there is a constant $C$ such that $||Tf||_{L^2(\R^n)} \leq C||f||_{L^2(\R^n)}$
for all Schwartz functions $f(x)$.

The condition concerning zeroes on lines parallel to the coordinates axes is needed for technical reasons in the proof. Note that
this condition holds whenever $b(x)$ is a polynomial, and it is not hard to see that it always holds in two variables, using the Weierstrass
Preparation Theorem for example. The author does not know
if it holds for all real-analytic functions, so it is included as an assumption in Theorem 2.1 (and in Theorem 2.2 below).

\noindent For $p \neq 2$, we have a weaker statement. To motivate the statement of the theorem, 
in $(4.2)$ and the line afterwards we will see that if each polynomial $b_F(x)$ is nonvanishing on $(\R - \{0\})^n$, then  there are constants $C_1$ and $C_2$ such that 
$$C_1b^*(x) < |b(x)| < C_2 b^*(x) \eqno (2.11)$$
 Hence in this situation $(2.4)$ becomes
$$|\partial_{y_l}K_{j_1,...,j_n}(y_1,...,y_n)| < C_1' {1 \over |y_l|} (b^*(y_1,...,y_n))^{-\delta_0}  \eqno (2.12)$$
For the $L^p$ theorem, we need bounds on derivatives of higher order in order to apply the Marcinkiewicz multiplier theorem. Hence
we assume that each $K_{j_1,...,j_n}(y_1,...,y_n)$ is a $C^{n+1}$ function and there is a constant $C$ such that for any multiindex $\alpha$
with $0 \leq |\alpha| \leq n + 1$ we have
$$|\partial^{\alpha} K_{j_1,...,j_n}(y_1,...,y_n)| \leq C {1 \over |y_1|^{\alpha_1}...|y_n|^{\alpha_n}}(b^*(y_1,...,y_n))^{-\delta_0} \eqno (2.13)$$
The condition $(2.13)$ is motivated by the fact that by iterating $(2.2a)$, the bounds $(2.13)$ hold for $|b(x)|^{-\delta_0}$ whenever each $b_F(x)$ is nonvanishing on $(\R - \{0\})^n$. 

\noindent Our $L^p$ theorem is as follows.

\noindent {\bf Theorem 2.2.} Suppose each polynomial $b_F(x)$ is nonvanishing on $(\R - \{0\})^n$. Suppose also that there is a $C_0 > 0$ and a neighborhood $U$ of the origin such that for each $l$, there is a set $Z \subset \R^{n-1}$ of measure zero such that the function $\partial_{x_l} b(x)$ has at most $C_0$ zeroes in $U$ 
on any line parallel to the $x_l$ coordinate axis whose projection onto the plane $x_l = 0$ is not in $Z$. Then there is an $R > 0$ such that if each $K_{j_1,...,j_n}(y)$ satisfies $(2.3),(2.13),(2.5)$ and is supported on $|y| < R$, then if $1 < p  < \infty$
there is a constant $C$  such that $||Tf||_{L^p(\R^n)} \leq C||f||_{L^p(\R^n)}$
for all Schwartz functions $f(x)$.

Going back to the examples preceding $(2.6)$, in the first example where $b(x) = x_1^{a_1}...\,x_n^{a_n}$, each kernel
$K(x) = (-1)^{sgn(x_1) + ... + sgn(x_n)} \phi(x_1^2,...,x_n^2) |b(x)|^{-\delta_0}$ will be covered by Theorems 2.1 and 2.2.
As for the second example where $b(x) = f(x_1^2,...,x_n^2)$, the maximum order of any zero of any $b_F(x)$ on $(\R - \{0\})^n$
is the same as the maximum order of any $f_F(x)$ on $(\R^+)^n$. So when this quantity is less than $d(b) = 2 d(f)$, $K(x)$
will fall under the conditions of Theorem 2.1. When each $f_F(x)$ is nonvanishing on $(\R^+)^n$, then $K(x)$ will fall under
the conditions of Theorem 2.2 as well.

\noindent {\bf 3. Proofs of theorems when the multiplicity is equal to one.}

\noindent Since Theorem 1.1 is a special case of Theorem 1.2, we prove Theorem 1.2. 

\noindent {\bf Proof of Theorem 1.2.}

Theorem 1.2 will follow if we can prove each that for each $i$ there is a constant $C$ such that $||T_if||_{L^2(\R^{m})} \leq C||f||_{L^2(\R^{m})}$ for all Schwartz functions $f(x)$. Here $T_i$  is as in $(1.15)$. Let $m_{ij}$ and $e_{ij}$ be the monomial exponents of $b \circ \beta_i (x)$ and $Jac_{\beta_i} (x)$ as before. If $i$ is such that $\delta_0 < {1 + e_{ij} \over m_{ij}}$ for each $j$, then the kernel of $T_i$ is
absolutely integrable and $L^2$ boundedness is immediate. Thus it suffices to consider only the $i$ for which there is some $l$ for which $\delta_0 = {1 + e_{il} \over m_{il}}$. Since we are assuming $b(x)$ has multiplicity one, there will only be one such $l$ for
each such $i$. Writing $h\circ \beta_i(y) = (h_1\circ \beta_i(y),...,h_m\circ \beta_i(y))$, $(1.15)$ can be rewritten as
$$T_i f(x) = \int_{\R^n} f(x - h \circ \beta_i(y))\,\alpha (x,\beta_i(y))\, (\rho_i \circ \beta_i(y))\, k_i(y)\, Jac_{\beta_i}(y)\,dy \eqno (3.1)$$
It is more convenient for our proofs that there be no nonzero $\lambda$ such that $\lambda \cdot (h\circ \beta_i(y))$ is a 
linear function of $y$. This can be accomplished as follows. If there is a nonzero $\lambda$ such that $\lambda \cdot (h\circ \beta_i(y))$ is the zero function, then on each hyperplane orthogonal to $\lambda$, the operator $T_i$ restricts to an operator of the same type as $T_i$ here, except the ambient space is of one lower dimension. Repeating as necessary, we may assume that 
$\lambda \cdot (h\circ \beta_i(y))$ is never the zero function for any nonzero $\lambda$. It is worth mentioning that we are
using the fact that $h \circ \beta_i(y)$ extends to a connected neighborhood of the support of $\rho_i \circ \beta_i(y)$ to ensure
we don't have different functions on different connected components to worry about.

One can further ensure that $\lambda \cdot (h\circ \beta_i(y))$ is never linear by letting $(y_1,...,y_n) = (z_1^3,...,z_n^3)$ with the corresponding 
change from $\beta_i(y_1,...,y_n)$ to $\tilde{\beta}_i(z) = \beta_i(z_1^3,...,z_n^3)$. The exponents
$m_{ij}$ and $e_{ij}$ can change, but $(1.6)-(1.8)$ and the other properties from resolution of singularities that we are using
 will still hold. Thus in the following, without loss of generality we will always assume that we are working in a situation where
 $\lambda \cdot (h\circ \beta_i(y))$ is not linear for any nonzero $\lambda$.

\noindent Let $k_{iL}(x) = \sum_{ (j_1,...,j_n): j_l < L \,\,for\,\,all\,\,  l} \,\,\,k_{i,j_1,...,j_n}(x)$ as in $(1.9)$, and let $T_{iL}$ be defined by
$$T_{iL}f(x) = \int_{\R^n} f(x - h \circ \beta_i(y)) \,\alpha (x,\beta_i(y))\, (\rho_i \circ \beta_i(y))\, k_{iL}(y)\, Jac_{\beta_i}(y)\,dy  \eqno (3.2)$$
In order to prove Theorem 1.2, it suffices to show that $T_{iL}$ is bounded on $L^2$ with bounds uniform in $L$. We can reduce to the case where $(\rho_i \circ \beta_i(y)) \alpha(x,\beta_i(y))$ is replaced by a function of $y$ (i.e. the operator is translation-invariant) through the following lemma.

\noindent {\bf Lemma 3.1.} Define  $U_{iL}$ by
$$U_{iL}f(x) = \int_{\R^n} f(x - h \circ \beta_i(y))\, k_{iL}(y) \,Jac_{\beta_i}(y)\,dy \eqno (3.3)$$
 If there is a constant $C$ depending on $b(x)$ (and the resolution of singularities procedure we are using on it), $h_1(x),...,h_m(x)$, and the constant $C_0$ of $(1.6)-(1.8)$, such that $||U_{iL}f||_p \leq C||f||_p$ 
for all Schwartz $f$ and all $L$, then there is a constant $C'$ such that $||T_{iL}f||_p \leq C'||f||_p$ for all Schwartz $f$ and all $L$.

\noindent {\bf Proof.} Let $\gamma(x,y)$ be the Schwartz function $(\rho_i \circ \beta_i(y)) \,\alpha(x,\beta_i(y))$. We use the Fourier inversion formula in the $x$ variable and write
$$\gamma(x,y) = \int_{ \R^{m}} \hat{\gamma}(t,y) \,e^{it_1 x_1 + ... + it_{m} x_{m}}\,dt \eqno (3.4)$$
Here $\hat{\gamma}(t,y)$ refers to the Fourier transform in the $x$ variable only. If $f$ and $g$ are Schwartz functions then $\int_{\R^n} T_{iL} f(x) g(x)\, dx $  is equal to
$$ \int_{\R^{m}}\bigg[\int_{\R^{m}}\bigg(\int_{\R^n} f (x - h \circ \beta_i(y)) \, \hat{\gamma}(t,y)\,k_{iL}(y)\,\,Jac_{\beta_i}(y) \,dy\bigg)$$
$$\times  (e^{it_1 x_1 + ... + it_{m} x_{m}}g(x))\,dx\bigg]\,dt \eqno (3.5)$$
Stated another way, let  $U_{iLt}$ denote the operator
$$U_{iLt}f(x) = \int_{\R^n} f (x - h \circ \beta_i(y))\, \hat{\gamma}(t,y)\,k_{iL}(y)\,Jac_{\beta_i}(y) \,dy \eqno (3.6)$$
Then $\int_{\R^n} T_{iL} f(x) g(x)\, dx $  is equal to
$$\int_{\R^{m}} \int_{\R^{m}} U_{iLt}f(x)\times (e^{it_1 x_1 + ... + it_{m} x_{m}}  g(x))\,dx \,dt \eqno (3.7)$$
Since $\hat{\gamma}(t,y)$ is Schwartz, under the assumptions of this lemma there is a constant $K$  such that 
$$|| U_{iLt}||_{L^p(\R^n) \rightarrow L^p(\R^n)}  \leq K{1 \over 1 + |t|^{m + 1}} \eqno (3.8)$$
Thus by $(3.7)$ and H\"older's inequality we have
$$\bigg|\int_{\R^n} T_{iL} f(x) g(x)\, dx \bigg| \leq K ||f||_p||g||_{p'}  \int_{\R^m}{1 \over 1 + |t|^{m + 1}}\, dt\eqno (3.9)$$
Thus the $T_{iL}$ are bounded on $L^p$  uniformly in $L$ and we are done with the proof of Lemma 3.1.

We now proceed to proving uniform bounds on the $U_{iL}$. Taking Fourier transforms, we get
$$\widehat{U_{iL}}(\lambda) = 
\hat{f}(\lambda) \int_{\R^n} e^{ i\lambda \cdot (h \circ \beta_i(y))}\, k_{iL}(y) \,Jac_{\beta_i}(y)\,dy\eqno (3.10)$$
Hence in order to prove Theorem 1.2, it suffices to show that there is a constant $C$ such that $|B_{iL}(\lambda)| \leq C$
for each $i$, $L$ and $\lambda$, where 
$$B_{iL}(\lambda)= \int_{\R^n} e^{ i\lambda \cdot (h \circ \beta_i(y))}\,k_{iL}(y)\, Jac_{\beta_i}(y)\,dy \eqno (3.11)$$
Without loss of generality, to simplify notation in the following we will assume that the $l$ for which  $\delta_0 = {1 + e_{il} \over m_{il}}$ is $l = 1$. Next, we write the factor  $ k_{iL}(y)  Jac_{\beta_i}(y)$ in $(3.11)$ as $\sum_{m < L} p_{i m L}(y)$, where
we add over the dyadic pieces in the $y_2,...,y_n$ variables to form $p_{imL}(y)$:
$$p_{i m L}(y) = \sum_{{ (j_2,...,j_n): j_l < L \,\,for\,\,all\,\,  l > 1}} k_{i,m,j_2,...,j_n}(y) \,Jac_{\beta_i}(y) \eqno (3.12)$$
Then $(1.6)$  implies the estimates
$$|p_{i m L}(y)| \leq C|y_1|^{-\delta_0m_{i1} + e_{i1}}...|y_n|^{-\delta_0m_{in} + e_{in}} \eqno (3.13)$$
Similarly, $(1.6)-(1.7)$ implies that for each $l$ we have
$$|\partial_{y_l}\, p_{i m L}(y)| \leq C {1 \over |y_l|} |y_1|^{-\delta_0m_{i1} + e_{i1}}...|y_n|^{-\delta_0m_{in} + e_{in}} \eqno (3.14)$$
 The cancellation condition $(1.8)$ gives
$$\bigg|\int_{\R} p_{i m L}(y_1,...,y_n) \, dy_1\bigg| < C 2^{-\epsilon_0 m} \eqno (3.15)$$
Note that $-\delta_0m_{i1} + e_{i1} = -1$  here, while the other exponents are all greater than $-1$. If one changes variables
$y_l = z_l^N$ in $(3.11)$ for some $l > 1$, instead of having a factor $|y_l|^{-\delta_0m_{il} + e_{il}}$ in $(3.13)-(3.14)$ one
has a factor of $|z_l|^{(-\delta_0m_{il} + e_{il})N}$. One also gains an additional factor of $N |z_l|^{N-1}$ from the Jacobian
of the coordinate change. Thus overall one has a factor of $|z_l|^{(-\delta_0m_{il} + e_{il} + 1)N - 1}$. 
Since $-\delta_0m_{il} + e_{il} > -1$, if $N$ is large enough this factor will be bounded by just $|z_l|$ and thus one can remove
the $z_l$ variable from $(3.13)-(3.14)$. Stated another way, if one changes $y_l = z_l^N$ for all $l > 1$ (making $N$ odd to ensure
it's a one-to-one map), $(3.11)$ becomes
$$B_{iL}(\lambda)= \int_{\R^n} e^{  i\lambda \cdot  (h \circ \beta_i)(z_1,z_2^N,...,z_n^N)}\,q_{iL}(z) \,dz \eqno (3.16)$$
Here $q_{iL}(z) = \sum_{m < L} q_{i m L}(z)$, where $q_{imL}(z)$ satisfies
$$|q_{i m L}(z)| \leq C|z_1|^{-1} \eqno (3.17)$$
$$|\partial_{z_1}\, q_{i m L}(z)| \leq C |z_1|^{-2},{\hskip 0.5in} \forall l > 1 \,\,\,\,|\partial_{z_l}\, q_{i m L}(z)| \leq C |z_1|^{-1}  \eqno (3.18)$$
$$\bigg|\int_{\R} q_{i m L}(z_1,...,z_n) \, dz_1\bigg| \leq C2^{-\epsilon_0 m} \eqno (3.19)$$
Letting  $f_i(z) =  h \circ \beta_i(z_1,z_2^N,...,z_n^N)$,  $(3.16)$ becomes
$$B_{iL}(\lambda)= \int_{\R^n} e^{i\lambda \cdot f_i(z)} q_{iL}(z)\,dz \eqno (3.20)$$
In view of the discussion above $(3.2)$, we may assume that $\lambda \cdot f_i(z)$ is not linear for any nonzero $\lambda$. Therefore,
writing $\lambda = |\lambda| \omega$ for $\omega \in S^{m-1}$, by a compactness argument on $S^{m-1} \times supp(q_{iL})$,
we may restrict consideration to $\omega$ to a small neighborhood $N$ in $S^{m-1}$ and replace $q_{iL}(x)$ by
$\sigma(x) q_{iL}(x)$ for a function $\sigma(x)$ supported on a ball $B(x_0,r_0)$ on which there is an $\epsilon > 0$ and a single directional
derivative $\partial_v$ such that  $|\partial_v^{\alpha} (\omega \cdot f_i)(z)| > \epsilon$ on $B(x_0,r_0)$ for some $\alpha \geq 2$. We can do
this in such a way that $v$ has a positive $x_1$ component. Thus if $\bar{x}$ denotes the first component of $x_0$, for $\omega \in N$ we are attempting to bound $\int_{\R^{n-1}}D_{iL}(\lambda,z_2,...,z_n)\,dz_2...dz_n$, where
$$D_{iL}(\lambda,z_2,...,z_n)= \int_{\R} e^{i|\lambda| (\omega \cdot f_i)((\bar{x},z_2,...,z_n) + tv)} \sigma((\bar{x},z_2,...,z_n) + tv)$$
$$\times  q_{iL}((\bar{x},z_2,...,z_n) + tv)\,dt  \eqno (3.21a)$$
Thus to prove Theorem 1.2, it suffices for our purposes to bound $D_{iL}(\lambda,z_2,...,z_n)$ uniformly in $L,\lambda,z_2,...,z_n$ for $\omega \in N$. For this, it suffices to bound $\tilde{D}_{iL}(\lambda,\tilde{\lambda}, z_2,...,z_n)$ uniformly in $L,\lambda, \tilde{\lambda},z_2,...,z_n$ for $\omega \in N$, where 
$$\tilde{D}_{iL}(\lambda,\tilde{\lambda}, z_2,...,z_n) = \int_{\R} e^{i|\lambda| (\omega \cdot f_i)((\bar{x},z_2,...,z_n) + tv) + i\tilde{\lambda}t} \sigma((\bar{x},z_2,...,z_n) + tv)$$
$$\times q_{iL}((\bar{x},z_2,...,z_n) + tv)\,dt  \eqno (3.21b)$$

But similarly to $(3.10)-(3.11)$, such  uniform bounds for $\tilde{D}_{iL}(\lambda,\tilde{\lambda}, z_2,...,z_n)$ 
follows from uniform boundedness on
$L^2$ of the singular Radon transforms along curves in $\R^2$ of the form
$$U_{i\,L\,\omega\, z_2...z_n} f(x_1,x_2) = \int_{\R} f\big(x_1 -t, x_2 - (\omega \cdot f_i)((\bar{x},z_2,...,z_n) + tv)\big)$$
$$\times \sigma((\bar{x},z_2,...,z_n) + tv)\,q_{iL}((\bar{x},z_2,...,z_n) + tv)\,dt \eqno (3.22)$$
Note that 
$$|q_{imL}((\bar{x},z_2,...,z_n) + tv)- q_{imL}((\bar{x},z_2,...,z_n) + t(v_1,0,...,0))| \leq C|t|\max_{l > 1} \sup_z | \partial_{z_l}q_{imL}(z)|\eqno (3.23)$$
By $(3.18)$, since $|t| <  C 2^{-m}$  and $|z_1| \sim 2^{-m}$ we have
$$|q_{imL}((\bar{x},z_2,...,z_n) + tv) - q_{imL}\big((\bar{x},z_2,...,z_n) + t(v_1,0...,0)\big)| \leq C' \eqno (3.24)$$
Hence by the cancellation condition $(3.19)$ one has
$$\int_{\R} q_{imL}((\bar{x},z_2,...,z_n) + t v)\,dt < C'' 2^{-\epsilon_0 m} \eqno (3.25)$$ 
Since $\sigma((\bar{x},z_2,...,z_n)) +t  v) =  \sigma(\bar{x},z_2,...,z_n) + O(|t|)$, using $(3.17)$ one also has
$$\int_{\R} \sigma((\bar{x},z_2,...,z_n) +t  v)\,q_{imL}((\bar{x},z_2,...,z_n) + t v)\,dt < C''' 2^{-\epsilon_0 m} \eqno (3.26)$$
In other words, we have a cancellation condition in $(3.22)$ derived from $(3.19)$. The constant $C'''$ in $(3.26)$ depends on
$b(x), h_1(x),...,h_m(x)$ and the constant $C$ of $(3.17)-(3.19)$, which in turn depends on $b(x)$ and the constant $C_0$ of $(1.6)-(1.8)$.

The arguments of [G3] provide $L^2$ bounds for the operators $U_{i\,L\,\omega\, z_2...z_n}$ under the assumptions
$(3.17)-(3.19)$ and a  lower bound on $|\partial_v^{\alpha} (\omega \cdot f_i)(z)|$ that are uniform in $L$, $\omega$,
$z_2,...,z_n$ for $\omega \in N$. (A slightly stronger cancellation condition is assumed but $(3.19)$ suffices). This is because the bounds obtained
in [G3] are at least as strong as the bounds obtained when the convolution is over the curve $(t, t^{\alpha})$, in which case the bounds
can be expressed in terms of the constant $C$ of $(3.17)-(3.19)$, the constant $C''$ of $(3.26)$, and the function $h \circ \beta_i(x)$.  For the
ball $B(x_0,r_0)$ on which $\sigma(x)$ is supported, how small $r_0$ needs to be for the uniform bounds to hold will also be uniform
in the various parameters but may be smaller than the $r_0$ we originally selected. However, this can be corrected by writing $\sigma(x)$ as a finite sum of bump functions with smaller support if needed. This completes the proof of Theorem 1.2.

\noindent {\bf 4. Proof of theorems when the multiplicity is greater than one.}

We start with some facts from [G1] - [G2] which will help us understand the distribution function of $b(x)$ and related 
properties of integrals of $|b(x)|^{-\delta_0}$. The constructions in [G1] are slightly better for our purposes so we bring our
attention to them. By Lemmas 3.2' - 3.5 of [G1], if $U$ is a sufficiently small neighborhood of the origin, up to a
set of measure zero one may write $U = \cup_{i = 1}^N U_i$ as a finite union of open sets such that the following hold. 
Each
$U_i$ is contained in one of the $2^n$ octants determined by coordinate hyperplanes. For each $i$, there is 
some integer $1 \leq k_i \leq n$ and 
a function $\gamma_i: \R^n \rightarrow \R^n$ such that each component of $\gamma_i(x)$ is plus or minus a monomial and
 $\gamma_i^{-1}(U_i)$ satisfies the following. If $k_i < n$, then there are cubes $(0,\eta_i)^{k_i}$ and $(0,\eta_i')^{k_i}$ with $\eta_i' > \eta_i$ and
bounded open sets $O_i \subset O_i'$ whose closures are a subset of $\{x \in \R^{n - k_i}: x_l > 0$ for all $l\}$, such that 
$$(0,\eta_i)^{k_i} \times O_i \subset \gamma_i^{-1}(U_i) \subset (0,\eta_i')^{k_i} \times O_i' \eqno (4.1a)$$
If $k_i = n$, then there are cubes  $(0,\eta_i)^{k_i}$ and $(0,\eta_i')^{k_i}$ with $\eta_i' > \eta_i$ such that 
$$(0,\eta_i)^{k_i} \subset \gamma_i^{-1}(U_i) \subset (0,\eta_i')^{k_i} \eqno (4.1b)$$
In either case, there is a monomial $m_i(x_1,...,x_{k_i})$ and constants $C_i, C_i'$ such that on $\gamma_i^{-1}(U_i)$ one has
$$ C_i m_i(x_1,...,x_{k_i}) < |b^* \circ \gamma_i (x)| < C_i' m_i(x_1,...,x_{k_i}) \eqno (4.2)$$
If $k_i = n$, then $(4.2)$ holds with $b^*\circ \gamma_i(x)$ replaced by $b \circ \gamma_i(x)$. If $k_i < n$, then
on $\gamma_i^{-1}(U)$  the function $b \circ \gamma_i(x)$ can be expressed as $m_i(x_1,...,x_{k_i})g_i(x_1,...,x_n)$ where
$m_i(x_1,...,x_{k_i})$ is a monomial and where $g_i(x_1,...,x_n)$ satisfies the following.
One may write $\gamma_i^{-1}(U_i) = \cup_{j = 1}^{M_i} V_{ij}$ such that for each $i$ and $j$ there is an $\epsilon > 0$,
a compact face $F_i$ of $N(b)$, and a directional
derivative $ \partial_{v_{ij}}$ in the last $n - k_i$ variables, such that $|\partial_{v_{ij}}^{a_{ij}}g_i(x_1,...,x_n)| > \epsilon$ on $V_{ij}$ for some $a_{ij} \geq 0$ which is at most the maximum order of any zero of $b_{F_i}(x_1,...,x_n)$ on $(\R - \{0\})^n$. 
When $a_{ij} = 0$, we interpret $\partial_{v_{ij}}^{a_{ij}}g_i(x)$ to just mean $g_i(x)$.

For a given $\epsilon > 0$, the following lemma  explicitly bounds  the measure of the portion of a dyadic rectangle where 
$|b(x)/b^*(x)| < \epsilon$  in terms of the maximum order of the zeroes of the $b_F(x)$ on $(\R - \{0\})^n$.

\noindent {\bf Lemma 4.1.} Suppose $p > 0$ is an integer such that the zeroes of each $b_F(x)$ on $(\R - \{0\})^n$ are all of order at most $p$. Then there is a neighborhood $U$ of the origin and a constant  $C > 0$ such that if $R \subset U$ is a set of the form $ \{x \in \R^n: 2^{-j_l} < |x_l| < 2^{-j_l+1}\}$ for integers $j_l$, then 
$$|\{x \in R: |b(x)/b^*(x)| < \epsilon\}| < C \epsilon^{1 \over p} |R| \eqno (4.3)$$
\noindent {\bf Proof.} It suffices to show for each $i$ an estimate of the form $|\{x \in R \cap U_i: |b(x)/b^*(x)| < \epsilon\}| < C \epsilon^{1 \over p} |R|$. Since the components of $\gamma_i(x)$ are all monomials, the absolute value of the Jacobian of
$\gamma_i(x)$ is of the form $c_i x_1^{e_{i1}}...\,x_n^{e_{in}}$ for some integers $e_{i1},...,e_{in}$ and some $c_i > 0$.
Viewing $|\{x \in R \cap U_i: |b(x)/b^*(x)| < \epsilon\}| $ as the integral of $1$ over $\{x \in R \cap U_i: |b(x)/b^*(x)| < \epsilon\}$
and changing coordinates via $\gamma_i(x)$, one obtains 
$$|\{ R \cap U_i : |b(x)/b^*(x)| < \epsilon\}| = \int_{\{x \in \gamma_i^{-1}(R) \cap  \gamma_i^{-1}(U_i):\,|b \circ \gamma_i (x)/b^* \circ \gamma_i(x)| < \epsilon \}} c_i x_1^{e_{i1}}...\,x_n^{e_{in}} \, dx \eqno (4.4)$$
Note that by $(4.2)$ and the following paragraph, one has $|b \circ \gamma_i (x)/b^* \circ \gamma_i(x)| > C' g_i(x)$ for some constant $C'$ (We can include the $k_i = n$ situation here by defining $g_i(x) = 1$). Thus in order to bound $(4.4)$ by an expression of the form $C\epsilon^{1 \over p} |R|$,
it suffices to show the following estimate of the following form for each $i$ and $j$.
$$ \int_{\{x \in \gamma_i^{-1}(R) \cap V_{ij}:\, g_i(x) < \epsilon \}} x_1^{e_{i1}}...\,x_n^{e_{in}} \, dx  < C\epsilon^{1 \over p} |R|\eqno (4.5)$$
If the multiindex $a_{ij}$ in the paragraph after $(4.2)$ is zero, then $g_i(x)$ is bounded below, and thus $(4.5)$ reduces to
showing that $ \int_{ \gamma_i^{-1}(R) \cap V_{ij}} x_1^{e_{i1}}...\,x_n^{e_{in}} \, dx $ is bounded by a constant times $|R|$, which follows immediately from changing back into the original coordinates using $\gamma_i$. Thus it suffices to assume $a_{ij} 
\geq 1$. Note that this only occurs if $k_i < n$. Since the final $n - k_i$ variables are bounded below on $V_{ij}$, it suffices to
prove a bound 
$$ \int_{\{x \in \gamma_i^{-1}(R) \cap V_{ij}:\, g_i(x) < \epsilon \}} x_1^{e_{i1}}...\,x_{k_i}^{e_{ik_i}} \, dx  < C\epsilon^{1 \over p} |R|\eqno (4.6)$$
We now integrate the left-hand side of $(4.6)$ starting with the $v_{ij}$ direction. Since
$a_{ij} \leq p$, by the
measure version of the Van der Corput lemma (see [C] for details), the integral in the $v_{ij}$ direction is at most $Cx_1^{e_{i1}}...\,x_{k_i}^{e_{ik_i}}\epsilon^{1 \over p}$. If we next perform the integration in the  remaining $n - k_i - 1$
directions of last $n - k_i$ 
variables (if any exist), then if $\pi_i$ denotes the projection on $\R^n$ onto the first
$k_i$ variables, we obtain
$$ \int_{\{x \in \gamma_i^{-1}(R) \cap V_{ij}:\, g_i(x) < \epsilon \}} x_1^{e_{i1}}...\,x_{k_i}^{e_{ik_i}} \, dx \leq
 C\epsilon^{1 \over p} \int_{\pi_i(\gamma_i^{-1}(R) \cap V_{ij})}  x_1^{e_{i1}}...\,x_{k_i}^{e_{ik_i}}\, dx_1...dx_{k_i} \eqno (4.7)$$
$$ = C\epsilon^{1 \over p} \int_{\pi_i(\gamma_i^{-1}(R) \cap V_{ij}) \times [1,2]^{n - k_i}} x_1^{e_{i1}}...\,x_{k_i}^{e_{ik_i}}\, dx_1...dx_{k_i} \eqno (4.8a)$$
$$ = C'\epsilon^{1 \over p} \int_{\pi_i(\gamma_i^{-1}(R) \cap V_{ij}) \times [1,2]^{n - k_i}} x_1^{e_{i1}}...\,x_{n}^{e_{in}}\, dx_1...dx_n \eqno (4.8b)$$
Because the last $n- k_i$ coordinates of the points in $U_i$ are bounded above and below away from zero, there is a constant
$C_0 > 1$ such that if $(x_1,...,x_n) \in \pi_i(\gamma_i^{-1}(R) \cap V_{ij}) \times [1,2]^{n - k_i}$ then there is a point $(y_1,...,y_n) \in \gamma_i^{-1}(R) \cap V_{ij}$ such that ${1 \over C_0} < {y_l \over x_l} < C_0$ for each $l$. This property
is preserved under monomial maps (perhaps with a different constant $C_1$), so the image of  $\pi_i(\gamma_i^{-1}(R) \cap V_{ij}) \times [1,2]^{n - k_i}$ under $\gamma_i$ is a subset of a corresponding dilate of $\gamma_i(\gamma_i^{-1}(R) \cap V_{ij})$, which
in turn is a subset of the dilate of $R$. Denote this dilate by $R^*$. Changing coordinates in $(4.8)$  back into the original coordinates via $\gamma_i$, we see that
$$C\epsilon^{1 \over p} \int_{\pi_i(\gamma_i^{-1}(R) \cap V_{ij}) \times [1,2]^{n - k_i}} x_1^{e_{i1}}...\,x_{n}^{e_{in}}\, dx_1...dx_n \leq C'' \epsilon^{1 \over p} \int_{R^*} 1 \,dx$$
$$= C'''\epsilon^{1 \over p}|R| \eqno (4.9)$$
This is the desired estimate $(4.6)$ and we are done.

\noindent We also will make use of the following result.

\noindent {\bf Lemma 4.2.}  Suppose the zeroes of each $b_F(x)$ on $(\R - \{0\})^n$ are all of order less than $d(b)$. Then there is a neighborhood $U$ of the origin and constants $C, \eta > 0$ such that if  $\epsilon > 0$ and $R \subset U$ is a set of the form $ \{x \in \R^n: 2^{-j_l} < |x_l| < 2^{-j_l+1}\}$, then $\int_{\{x \in R:\, |b(x)| < \epsilon |b^*(x)|\}}|b(x)|^{-\delta_0} < 
C\epsilon^{\eta}$. In particular, since there is a constant $C'$ such that $|b(x)| \leq C' b^*(x)$ on any such $R$, there is a constant $C''$ such that
$\int_R |b(x)|^{-\delta_0} < C''$ for such $R \subset U$.

\noindent {\bf Proof.} Since the terms of $b^*(x)$ are absolute values of monomials, there is a constant $c > 1 $ and an $x_0 \in R$ such that $c\, b^*(x_0) \geq b^*(x) \geq {1 \over c} b^*(x_0)$ on $R$. Hence it suffices to prove an
estimate of the form $\int_{\{x \in R:\,|b(x)| < \epsilon |b^*(x_0)|\}}|b(x)|^{-\delta_0} < C\epsilon^{\eta}$.
By the relation between $L^p$ norms and distribution functions, applied to ${1 \over |b(x)|}$, one has
$$\int_{\{x \in R:\,|b(x)| < \epsilon |b^*(x_0)|\}} |b(x)|^{-\delta_0} = \delta_0 \int_0^{\infty} t^{\delta_0 - 1}\big|\{x \in R: |b(x)| < \min( \epsilon |b^*(x_0)|, {1 \over t})\}\big|\,dt 
\eqno (4.10)$$
It is natural to break up $(4.10)$ into two pieces, the first where $t <{1 \over \epsilon b^*(x_0)}$ and the second where $t \geq {1 \over \epsilon b^*(x_0)}$. Then the right-hand side of $(4.10)$ becomes
$$\delta_0 \int_0^{ 1 \over  \epsilon b^*(x_0)} t^{\delta_0 - 1}|{\{ x \in R: \,|b(x)| < \epsilon |b^*(x_0)|\}}|\,dt +  \delta_0 \int_{ 1 \over  \epsilon b^*(x_0)}^{\infty} t^{\delta_0 - 1}\big|\{x \in R: |b(x)| < {1 \over t}\}\big|\,dt \eqno (4.11)$$
Performing the first integral in the first term of $(4.11)$ results in 
$$\epsilon^{-\delta_0} |b^*(x_0)|^{-\delta_0}|{\{ x \in R: \,|b(x)| < \epsilon |b^*(x_0)|\}}| \eqno (4.12)$$
By Lemma 4.1, $(4.12)$ is bounded by $C|b^*(x_0)|^{-\delta_0}\epsilon^{{1 \over p }- \delta_0}|R|$ for some $p < d(b) = {1 \over \delta_0}$. Hence we have
$$\int_{\{x \in R:\,|b(x)| < \epsilon |b^*(x_0)|\}} |b(x)|^{-\delta_0} \leq C|b^*(x_0)|^{-\delta_0}\epsilon^{{1 \over p }- \delta_0}|R|$$
$$+  \delta_0 \int_{ 1 \over \epsilon b^*(x_0)}^{\infty} t^{\delta_0 - 1}
\big|\{x \in R: |b(x)| < {1 \over t}\}\big|\,dt \eqno (4.13)$$
Note that $\{x \in R: |b(x)| < {1 \over t}\} \subset \{x \in R: |b(x)| < {c_0 \over tb^*(x_0)}b^*(x)\}$, so by Lemma 4.1 
for some constant $C_0$ we have
$$\int_{\{x \in R:\,|b(x)| < \epsilon |b^*(x_0)|\}} |b(x)|^{-\delta_0} \leq C|b^*(x_0)|^{-\delta_0}\epsilon^{{1 \over p }- \delta_0}|R|+  C_0\int_{ 1 \over  \epsilon b^*(x_0)}^{\infty} t^{\delta_0 - 1}
\bigg({1 \over t b^*(x_0)}\bigg)^{1 \over p}|R| \,dt \eqno (4.14)$$
Note that the exponent $\delta_0 - 1 - {1 \over p}$ is less than $-1$ since $p <  d(b) = {1 \over \delta_0}$. Hence integrating
the second term on the right of $(4.14)$  leads to the following for some constant $C_1$.
$$\int_{\{x \in R:\,|b(x)| < \epsilon |b^*(x_0)|\}} |b(x)|^{-\delta_0} \leq C |b^*(x_0)|^{-\delta_0}\epsilon^{{1 \over p }- \delta_0}|R| +  C_1 |b^*(x_0)|^{-\delta_0}\epsilon^{{1 \over p}- \delta_0}|R|\eqno (4.15)$$
Since $|R| \sim |x_1...\,x_n|$ for any $(x_1,...,x_n) \in R$, in order to prove Lemma 4.2 with $\eta = {1 \over p} - \delta_0$, it suffices to show that there is a 
constant $C_2$ such that for any $x$ we have
$$|x_1...\,x_n|( b^*(x))^{-\delta_0}  \leq C_2 \eqno (4.16)$$
Since $\delta_0 = {1 \over d(b)}$ in the case at hand, $(4.16)$ is equivalent to the statement that
$$b^*(x) \geq C_3 |x_1...\,x_n|^{d(b)} \eqno (4.17)$$
Since $(d(b),...,d(b))$ is on the Newton polyhedron $N(b)$, there are nonnegative $\alpha_i$ with $\alpha_1 + ... + \alpha_k =
 1$ such that each component of $(d(b),...,d(b))$ is greater than or equal to that of $\alpha_1 v_1 + ... + \alpha_k v_k$ for 
some vertices  $v_1,...,v_k$ of $N(b)$.  Hence $|x_1...\,x_n|^{d(b)} \leq  |x^{v_1}|^{\alpha_1}...|x^{v_k}|^{\alpha_k}$. So by the generalized AM-GM inequality one has  $|x_1...\,x_n|^{d(b)}\leq \sum_{i = 1}^k
\alpha_i |x^{v_i}| \leq b^*(x)$ as needed. This completes the proof of Lemma 4.2.

Similar to the multiplicity one case, in order to show 
$||Tf||_p \leq C||f||_p$ for all Schwartz $f$ for a given $1 < p < \infty$, it suffices to show the that if $K(y)$ is supported on a sufficiently small neighborhood of the origin there is a constant $C$ such that
$||T_Lf||_p \leq C||f||_p$ for all Schwartz $f$ and each $L$, where $T_L f(x) = \int_{\R^n} f(x - y) \,\alpha(x,y)\,K_L(y)\,dy$.
Here $K_L(y) = \sum_{j_l < L\,\,for\,\, all\,\, l} K_{j_1,...,j_n}(y)$ as in $(2.6)$. As in Lemma 3.1 for the multiplicity one case, we may also replace $\alpha(x,y)$ by just $1$. Thus we
focus our attention on $U_L$ given by
$$U_L f(x) = \int_{\R^n} f(x - y) K_L(y)\,dy \eqno (4.18)$$
Our goal will be to prove $U_L$ is bounded on $L^p$ with a norm independent of $L$ under the hypotheses of Theorem 2.1 or
2.4.
The next two lemmas provide bounds on the $|\widehat{K_{j_1,...,j_n}}(\xi)|$ that allow us to prove such uniform bounds.

\noindent {\bf Lemma 4.3.} Under the assumptions of Theorem 2.1, there is a constant $ C > 0$ such that if $l$ is such that ${2^{-j_l}|\xi_l|} \leq 1$, then 
$$|\widehat{K_{j_1,...,j_n}}(\xi)| \leq C 2^{-j_l} |\xi_l|$$
\noindent {\bf Proof.} $\widehat{K_{j_1,...,j_n}}(\xi)$ is given by
$$\widehat{K_{j_1,...,j_n}}(\xi) = \int_{\R^n} K_{j_1,...,j_n}(x) e^{-i\xi_1 x_1 - ...- i\xi_n x_n}\, dx \eqno (4.19)$$
Since the integral of $K_{j_1,...,j_n}(x)$ in the $x_l$ variable is equal to zero by $(2.5)$, one can subtract $K_{j_1,...,j_n}(x) e^{\sum_{k \neq l}-i\xi_k x_k}$ from the integrand in $(4.19)$ without changing the integral, so we have
$$\widehat{K_{j_1,...,j_n}}(\xi) = \int_{\R^n} K_{j_1,...,j_n}(x) (e^{-i\xi_l x_l} - 1)e^{\sum_{k \neq l}-i\xi_k x_k}\, dx \eqno (4.20)$$
Since $|\xi_l x_l| \sim 2^{-j_l}|\xi_l| \leq C$ when $K_{j_1,...,j_n}(x) \neq 0$, in $(4.20)$ one has that  $(e^{-i\xi_l x_l} - 1)\leq C|x_l \xi_l| < C'2^{-j_l}|\xi_l|$ and we get
$$|\widehat{K_{j_1,...,j_n}}(\xi)| \leq C'2^{-j_l}|\xi_l|\int_{\R^n} |K_{j_1,...,j_n}(x)| \, dx \eqno (4.21)$$
Using Lemma 4.2 we obtain the desired estimate
$$|\widehat{K_{j_1,...,j_n}}(\xi)| \leq  C''2^{-j_l}  |\xi_l| \eqno (4.22)$$
\noindent {\bf Lemma 4.4.} Under the assumptions of Theorem 2.1, there are constants $\rho, C > 0$ such that if $l$
is such that ${2^{-j_l}|\xi_l|} \geq 1$, then 
$$|\widehat{K_{j_1,...,j_n}}(\xi)| \leq C {1 \over ({2^{-j_l}|\xi_l|)}^{\rho}}  \eqno (4.23)$$
\noindent {\bf Proof.} Let $\sigma_1(x)$ be a smooth increasing nonnegative function on $\R^+$ with $\sigma_1(x) = 1$ for $|x|  < 1$ and 
$\sigma_1(x) = 0$ for $|x|  >  2$. Let $\sigma_2(x) = 1 - \sigma_1(x)$. For a constant $\rho_0 > 0$ to be determined by our arguments, for any fixed $x_0$ in the dyadic rectangle corresponding to $(j_1,...,j_n)$ we write
$$\widehat{K_{j_1,...,j_n}}(\xi) = \int_{\R^n}\sigma_1\bigg((2^{-j_l}|\xi_l|)^{\rho_0}{|b(x)| \over b^*(x_0)}\bigg) K_{j_1,...,j_n}(x) e^{-i\xi_1 x_1 - ...- i\xi_n x_n}\, dx$$
$$+   \int_{\R^n}\sigma_2\bigg((2^{-j_l}|\xi_l|)^{\rho_0}{|b(x)| \over b^*(x_0)}\bigg) K_{j_1,...,j_n}(x) e^{-i\xi_1 x_1 - ...- i\xi_n x_n}\, dx \eqno (4.24)$$
The first term of $(4.24)$ is bounded by
$$ \int_{\{x: |b(x)| \leq 2 (2^{-j_l}|\xi_l|)^{-\rho_0}\,\,b^*(x_0)\}}| K_{j_1,...,j_n}(x)| \eqno (4.25)$$
Using $(2.3)$ and Lemma 4.1, we see that this term is at most $C (2^{-j_i}|\xi_i|)^{-{\rho_0 \over d}}$ for some $d$, which
gives the bound of the right-hand side of $(4.23)$.

Proceeding to the second term of $(4.24)$, we integrate by parts, integrating the $e^{-i\xi_1 x_1 - ...- i\xi_n x_n}$ factor
in the $x_l$ variable and differentiating the remaining factors. The resulting term is given by
 $${1 \over i\xi_l} \int_{\R^n}\partial_{x_l} \bigg[\sigma_2\bigg((2^{-j_l}|\xi_l|)^{\rho_0}{|b(x)| \over b^*(x_0)}\bigg) K_{j_1,...,j_n}(x) \bigg]e^{-i\xi_1 x_1 - ...- i\xi_n x_n}\, dx \eqno (4.26)$$
If the $x_l$ derivative in $(4.26)$ lands on the $K_{j_1,...,j_n}(x)$ factor, one obtains a term which by $(2.4)$ is bounded by
$$C_1 {1 \over |\xi_l|} \int_{\R^n}\sigma_2\bigg((2^{-j_l}|\xi_l|)^{\rho_0}{|b(x)| \over b^*(x_0)}\bigg){1 \over |x_l|} b^*(x)|b(x)|^{-1 - \delta_0}\, dx\eqno (4.27)$$
Due to the $\sigma_2$ factor in $(4.27)$, on the support of the integrand of $(4.27)$ we have $|b(x)| \geq
 (2^{-j_l}|\xi_l|)^{-\rho_0}b^*(x_0)$. Thus $(4.27)$ is bounded by
$$C_1 {1 \over |\xi_l|} \int_{R_{j_1,...,j_n}}\sigma_2\bigg((2^{-j_l}|\xi_l|)^{\rho_0}{|b(x)| \over b^*(x_0)}\bigg){1 \over |x_l|} b^*(x)(2^{-j_l}|\xi_l|)^{\rho_0(1 + \delta_0)}b^*(x_0)^{-1-\delta_0}\, dx\eqno (4.28)$$
Here $R_{j_1,...,j_n}$ denotes the (expanded) dyadic rectangle-like set on which $K_{j_1,...,j_n}$ is supported. Since $\sigma_2(t)
\leq 1$ for all $t$, $|x_l| \sim 2^{-j_l}$ on $R_{j_1,...,j_n}$, and $b^*(x)$ is within a constant factor of $b^*(x_0)$ on $R_{j_1,...,j_n}$, $(4.28)$ is bounded by
$$C_2 \bigg({1 \over 2^{-j} |\xi_l|}\bigg) (2^{-j_l}|\xi_l|)^{\rho_0(1 + \delta_0)}\int_{R_{j_1,...,j_n}}(b^*(x))^{-\delta_0}\,dx \eqno (4.29)$$
By Lemma 4.2 (which applies to negative powers of the smaller function $|b(x)|$), we see that the above is bounded by
$$C_3 \bigg({1 \over 2^{-j}|\xi_l|}\bigg) (2^{-j_l}|\xi_l|)^{\rho_0(1 + \delta_0)} \eqno (4.30)$$
Thus so long as $\rho_0$ is chosen so that $\rho_0(1 + \delta_0) < {1 \over 2}$ for example, this term of $(4.26)$ satisfies
the bounds needed in this lemma.

\noindent We now bound the term where the derivative in $(4.26)$ lands on the $\sigma_2((2^{-j_l}|\xi_l|)^{\rho_0}{|b(x)| \over b^*(x_0)})$ factor. Observe that
$$\partial_{x_l} \bigg(\sigma_2\bigg((2^{-j_l}|\xi_l|)^{\rho_0}{|b(x)| \over b^*(x_0)}\bigg)\bigg) = \pm\bigg( (2^{-j_l}|\xi_l|)^{\rho_0}{\partial_{x_l} b(x) \over b^*(x_0)}\bigg) \sigma_2' \bigg((2^{-j_l}|\xi_l|)^{\rho_0}{|b(x)| \over b^*(x_0)}\bigg) \eqno (4.31)$$
Since $|b(x)| \geq (2^{-j_l}|\xi_l|)^{-\rho_0}b^*(x_0)$ in the support of the $\sigma'$ factor, by $(2.3)$, on the support
of the integrand of this term of $(4.26)$ we have
$$|K_{j_1,...,j_n}(x)| \leq C_4 (2^{-j_l}|\xi_l|)^{\rho_0 \delta_0}|b^*(x_0)|^{-\delta_0} \eqno (4.32)$$
As a result, the absolute value of the term of $(4.26)$ in question is bounded by
$$C_4|\xi_l|^{-1} (2^{-j_l}|\xi_l|)^{\rho_0 \delta_0}|b^*(x_0)|^{-\delta_0} \int_{R_{j_1,...,j_n}}\bigg|\partial_{x_l} \bigg(\sigma_2\bigg((2^{-j_l}|\xi_l|)^{\rho_0}{|b(x)| \over b^*(x_0)}\bigg)\bigg) \bigg| \, dx \eqno (4.33)$$
We first integrate in the $x_l$ variable  in $(4.33)$. By the hypotheses of Theorem 2.1 concerning zeroes of $\partial_{x_l}b(x)$,  for any fixed value of the remaining $n-1$ variables (outside a set of measure zero) the domain of integration in the $x_l$ variable can be written as the union of boundedly many intervals on which $\partial_{x_l}(\sigma_2((2^{-j_l}|\xi_l|)^{\rho_0}{|b(x)| \over b^*(x_0)}))$ does not
change sign. Thus on each of these intervals this derivative integrates back to the function. Since
$\sigma_2$ is bounded this means the $x_l$ integrals in $(4.33)$ are uniformly bounded in the remaining variables. Thus
doing the $x_l$ integral first and then integrating over the remaining variables shows that $(4.33)$ is bounded by
$$C_5  |\xi_l|^{-1} (2^{-j_l}|\xi_l|)^{\rho_0 \delta_0}|b^*(x_0)|^{-\delta_0}2^{\sum_{i \neq l} -j_i}\eqno (4.34)$$
Since $b^*(x) \sim b^*(x_0)$ on $R_{j_1,...,j_n}$, $(4.34)$ is bounded by
$$C_6  |\xi_l|^{-1}(2^{-j_l}|\xi_l|)^{\rho_0 \delta_0}2^{j_l} \int_{R_{j_1,...,j_n}}|b^*(x)|^{-\delta_0}\,dx\eqno (4.35)$$
As in $(4.30)$, the integral in $(4.35)$ is uniformly bounded and we obtain the bound
$$C_7  (2^{-j_l}|\xi_l|)^{\rho_0 \delta_0}{1 \over  2^{-j_1}|\xi_l|} \eqno (4.36)$$
So as long as $\rho_0\delta_0 < 1$, we see from $(4.36)$ that the term of $(4.26)$ under consideration is also is bounded by the right-hand side of $(4.23)$. We have now shown that the first term of $(4.24)$ and each term of $(4.26)$ all are bounded by the right-hand side of $(4.23)$ and thus we are done with the proof of Lemma 4.4.

\noindent {\bf Proof of Theorem 2.1.}

We will prove $L^2$ boundedness of $U_L$ uniformly in $L$  by bounding the Fourier transform $\widehat{K_L}(\xi)$ uniformly in $L$ and $\xi$. Since $\widehat{K_L}(\xi) = \sum_{(j_1,...,j_n) \in \Z^n: {j_l < L{\rm\,\,for\,\, all\,\,} l}}\widehat{K_{j_1,...,j_n}}(\xi)$, we have the bound
$$|\widehat{K_L}(\xi)| \leq \sum_{(j_1,...,j_n) \in \Z^n: {j_l  < L{\rm\,\,for\,\, all\,\,} l}}|\widehat{K_{j_1,...,j_n}}(\xi)| \eqno 
(4.37)$$
We use the better of the two estimates from Lemmas 4.2 and 4.3 in each term of $(4.37)$ then add the result. Let $(k_1,...,k_n)$ be the vector of integers
such that for each $l$, $2^{k_l}$ is the nearest power of 2 to $|\lambda_l|$. For any $M$ the number of $(j_1,...,j_n)$ such that 
$\max_l |j_l - k_l| = M$ is bounded by $CM^{n-1}$, and for each such $(j_1,...,j_n)$ Lemma 4.3 or 4.4 gives a bound
$|\widehat{K_{j_1,...,j_n}}(\xi)| \leq C' 2^{-\rho_1 M}$ for some $\rho_1 > 0$. Hence in $(4.37)$ the sum over all terms with 
$\max_l |j_l - k_l| = M$ is bounded by $C'' M^{n-1}2^{-\rho_1 M}$. Adding over all $M$ gives a uniform bound and we are done.

\noindent {\bf Proof of Theorem 2.2.}

We will make use of the Marcinkiewicz multiplier theorem (see Theorem 6' on p. 109 of [S]), which implies that  $L^p$ bounds on
$U_L$ that are uniform in $L$ will follow if we can show that there is a constant $C$ such that for each multiindex $\alpha$ with 
$|\alpha| \leq n$ and each $L$ we have the estimate
$$|\xi_1^{\alpha_1}...\xi_n^{\alpha_n} \partial^{\alpha} \widehat{K_L}(\xi_1,...,\xi_n)| \leq C \eqno (4.38)$$
Returning to the $x$ variables, this will follow as in the proof of Theorem 2.1  if we can show that for each multiindex $\alpha$ with $0 \leq |\alpha| \leq n$ the kernel $\partial^{\alpha}( x^{\alpha}K_{j_1,...,j_n}(x))$ satisfies the conditions of Lemmas 4.3 and 4.4. But the fact that $(2.13)$ holds for $K_{j_1,...,j_n}(x)$ immediately implies that $(2.13)$ also holds for $\partial^{\alpha}( x^{\alpha}K_{j_1,...,j_n}(x))$. The cancellation condition $(2.5)$ also holds for $\partial^{\alpha}( x^{\alpha}K_{j_1,...,j_n}(x))$; if the $x_l$ variable is not represented in $\alpha$ it can be shown by multiplying $(2.5)$ through by  $x^{\alpha}$ and
then applying $\partial^{\alpha}$ under the integral sign, while if the $x_l$ variable is represented in $\alpha$, then the integral $(2.5)$
is zero simply because one is integrating the derivative of a compactly supported $C^1$ function. Hence each kernel $\partial^{\alpha}( x^{\alpha}K_{j_1,...,j_n}(x))$ satisfies Lemmas 4.3 and 4.4  and Theorem 2.2 follows.

\noindent {\bf 5. References.}

\noindent [AGV] V. Arnold, S. Gusein-Zade, A. Varchenko, {\it Singularities of differentiable maps},
Volume II, Birkhauser, Basel, 1988. \parskip = 3pt\baselineskip = 3pt

\noindent [C] M. Christ, {\it Hilbert transforms along curves. I. Nilpotent groups}, Annals of Mathematics
(2) {\bf 122} (1985), no.3, 575-596.

\noindent [FS] R. Fefferman, E. Stein, {\it Singular integrals on product spaces}, Adv. in Math. {\bf 45} (1982), no. 2, 117-143.

\noindent [G1] M. Greenblatt, {\it A constructive elementary method for local resolution of singularities}, submitted.

\noindent [G2] M. Greenblatt, {\it Oscillatory integral decay, sublevel set growth, and the Newton
polyhedron}, Math. Annalen {\bf 346} (2010), no. 4, 857-895.

\noindent [G3] M. Greenblatt, {\it A method for proving $L^p$ boundedness of singular Radon transforms in codimension one for
 $1 < p < \infty$}, Duke Math J. {\bf 108} (2001) 363-393.

\noindent [H1] H. Hironaka, {\it Resolution of singularities of an algebraic variety over a field of characteristic zero I}, 
 Ann. of Math. (2) {\bf 79} (1964), 109-203.

\noindent [H2] H. Hironaka, {\it Resolution of singularities of an algebraic variety over a field of characteristic zero II},  
Ann. of Math. (2) {\bf 79} (1964), 205-326. 

\noindent [N] A. Nagel, {\it 39th annual Spring Lecture Series: Multiparameter Geometry and Analysis}, University of Arkansas, 2014

\noindent [NS] A. Nagel, E. Stein, {\it On the product theory of singular integrals}, Rev. Mat. Iberoamericana {\bf 20} (2004),
 no. 2, 531-561.

\noindent [S] E. M. Stein, {\it Singular integrals and differentiability properties of functions}, Princeton Mathematical Series, No. 30, Princeton University Press, Princeton, N.J., 1970.

\noindent [V] A. N. Varchenko, {\it Newton polyhedra and estimates of oscillatory integrals}, Functional 
Anal. Appl. {\bf 18} (1976), no. 3, 175-196.

\line{}
\line{}

\noindent Department of Mathematics, Statistics, and Computer Science \hfill \break
\noindent University of Illinois at Chicago \hfill \break
\noindent 322 Science and Engineering Offices \hfill \break
\noindent 851 S. Morgan Street \hfill \break
\noindent Chicago, IL 60607-7045 \hfill \break

\end